\documentclass[12pt,a4paper]{article}

\usepackage{amsmath, amssymb}
\usepackage{graphicx}

\usepackage{color,verbatim}

\newcommand{\N}{\mathbb N}
\newcommand{\R}{\mathbb R}
\newcommand{\Z}{\mathbb Z}
\newcommand{\calF}{\mathcal{F}}

\newcommand{\bx}{\boldsymbol{x}}
\newcommand{\bk}{\boldsymbol{k}}
\newcommand{\br}{\boldsymbol{r}}

\title{Numerical solutions to an inverse problem for a non-linear Helmholtz equation}

\author{Q.~T.~Le Gia$^1$ \and H.~N.~Mhaskar$^2$}
\date{%
	$^1$ School of Mathematics and Statistics, \\UNSW, Sydney, Australia \\%
	$^2$ Institute of Mathematical Sciences, \\Claremont Graduate University, U.S.A \\[2ex]
    }

\begin{document}
\maketitle

\begin{abstract}
In this work, we construct numerical solutions to an inverse problem 
of a nonlinear Helmholtz equation defined in a spherical shell 
between two concentric spheres centred at the origin.
Assuming that the values of the forward problem are known at sufficiently 
many points, we would like to determine the form of the non-linear term on 
the right-hand side of the equation via its Chebyshev coefficients. 
\end{abstract}
\tableofcontents

\section{Introduction}\label{sec:intro}

The nonlinear Helmholtz equation (NLH) models the propagation of electromagnetic 
waves in Kerr media, and describes a range of important phenomena in nonlinear 
optics and in other areas \cite{FibTsy01,FibTsy05,BarFibTsy07}.
In this paper, we consider forward and inverse problems regarding the following nonlinear Helmholtz equation in $\R^3$:
\begin{equation}\label{equ:Helm}
 \Delta U(\bx) + k^2  \nu(\bx) U(\bx) =  -\epsilon(\bx) F(|U(\bx)|^2) U(\bx),
\quad \bx \in \Omega \subset \R^3. 
\end{equation}
where $\bx=(x_1,x_2,x_3)$ are the spatial coordinates, $U = U(\bx)$ denotes
the scalar electric field, $|\cdot|$ denotes the Euclidean norm,
$\Delta = \partial^2_{x_1} + \partial^2_{x_2} + \partial^2_{x_3}$ is
the Laplacian operator, $\nu(\bx)$ and $\epsilon(\bx)$ are some functions. 

For simplicity, we consider the case where $\Omega$ is a spherical shell 
between two concentric spheres of radii $R_0$ and $R_1$ centred at the origin, 
that is
\[
  \Omega := \{ \bx \in \R^3:    R_0 \le |\bx| \le R_1 \}.
\]
We also assume that $\nu$ and $\epsilon$ are radially symmetric and the $U$
satisfies the axially symmetric boundary conditions
\begin{equation}\label{equ:boundary}
U|_{r=R_0} = H(t), \qquad 
\left. \frac{ \partial U } {\partial r } \right|_{r=R_0} = G(t), \quad -1 \le t \le 1,
\end{equation}
where $r = |\bx|$ and $t = \cos \theta$, with $\theta$ being the polar angle
measured from the north pole. 
The solution $U$ is then axially symmetric as well.
Equation~\eqref{equ:Helm} takes the form
\begin{equation}
\label{eq:helmradial}
\frac{1}{r}\frac{\partial^2}{\partial r^2}(rU(r,t))+\frac{1}{r^2}\Delta^{(L)}U(r,t)+ k^2  \nu(r) U(r,t) =  -\epsilon(r) F(|U(r,t))|^2) U(r,t),
\end{equation}
where $\Delta^{(L)}$ is the Legendre differential operator defined below in \eqref{eq:legendreeqn}.

In the forward problem, $U$ is unknown, $F$ is non-linear, e.g.  $F(|U|^2) = |U|^{2p}$ with
some integer $p$, or $F=\sin(|U|^2)$, and we find an approximation of $U$.
In the inverse problem, the values of the solution $U(\bx_q)$, for $q=1,\ldots,Q$ are known, and the problem is to approximate the unknown nonlinear function $F$. 

The paper is organized as follows. In Section~\ref{sec:forward} we introduce a spectral method
for the forward problem and a fast algorithm to evaluate the non-linear term. In Section~\ref{sec:inverse}
we describe an algorithm for the inverse problem to identify the nonlinearity of the $F$ via its
Chebyshev coefficients. The paper is concluded with some numerical experiments described
in Section~\ref{sec:numerics}.

\section{Background}\label{sec:back}

The Legendre polynomial $P_\ell$ is a polynomial of degree $\ell$ with leading coefficients. 
We have the orthogonality relation
\begin{equation}
\label{eq:legendreortho}
\int_{-1}^1P_\ell(t)P_{\ell'}(t)dt=\frac{2\ell+1}{2}\delta_{\ell,\ell'}.
\end{equation}
The polynomials $P_\ell$ satisfy
\begin{equation}
\label{eq:legendreeqn}
\Delta^{(L)}P_\ell(t)=(1-t^2)P_\ell''(t)-2tP_\ell'(t)=-\ell(\ell+1)P_\ell(t)=-\lambda_\ell P_\ell(t).
\end{equation}
The Fourier-Legendre coefficients of an integrable function $g: [-1,1]\to \R$ are defined by
\begin{equation}
\label{eq:legendrefour}
\hat{g}(\ell)=\int_{-1}^1 g(t)P_\ell(t)dt, \qquad \ell\in \Z_+,
\end{equation}
To compute the Fourier-Legendre coefficients of a product of two functions, we define
\begin{equation}
\label{eq:legendreconn}
\Gamma(L;\ell,\ell')=\frac{2L+1}{2}\int_{-1}^1 P_L(t)P_\ell(t)P_{\ell'}(t)dt.
\end{equation}
It is known that $0\le \Gamma(L;\ell,\ell')\le 1$ and $\sum_{L=0}^{\ell+\ell'} \Gamma(L;\ell,\ell')=1$,
see \cite[Chapter~5]{askeyspecial}. 
Obviously, the following formal equation holds:
\begin{equation}
\label{eq:jacobiconv}
\widehat{g_1g_2}(L)=\sum_{\ell,\ell'}\Gamma(L;\ell,\ell')\hat{g_1}(\ell)\hat{g_2}(\ell').
\end{equation}
In terms of the sequences of Fourier-Legendre coefficients, we denote
\begin{equation}
\label{eq:jacobiconvfour}
\left(\hat{g_1}\star \hat{g_2}\right)(L)=\sum_{\ell,\ell'}\Gamma(L;\ell,\ell')\hat{g_1}(\ell)\hat{g_2}(\ell').
\end{equation}

\section{Spectral method for the forward problem}\label{sec:forward}
In this section,  we discuss how to construct a numerical solution to \eqref{equ:Helm}. 
For this purpose, we first establish some notation.

The spectral method for the forward problem is to find an approximation $U_N$ defined by
\begin{equation}\label{AxSuN}
   U_N(r, t) = \sum_{\ell=0}^N  u_{\ell}(r) P_{\ell}(t),\qquad u_\ell(r)=\widehat{u(r,\cdot)}(\ell),
\end{equation}
and find the coefficients $u_\ell(r)$ so that $U_N$ satisfies \eqref{eq:helmradial}.
 By substituting $U_N$ into \eqref{eq:helmradial},  we deduce using \eqref{eq:legendreeqn} that
\begin{equation}\label{ode:uell}
\frac{2}{2\ell+1} \left( \frac{1}{r} \frac{\partial^2}{\partial^2 r} (r u_{\ell} ) 
 - \frac{\lambda_\ell}{r^2} u_{\ell}  + k^2 \nu(r)  u_{\ell} \right) = 
-\epsilon(r)  \calF_{\ell},
\end{equation}
where
\begin{equation}\label{eq:Fell}
 \calF_{\ell} = \calF_{\ell}(r)= \int_{-1}^{1} U_N(r,t) F(|U_N(r,t)|^2) P_{\ell}(t) dt .
\end{equation}
Equivalently,
\begin{equation}\label{ODE_Fell}
\frac{ \partial^2}{\partial r^2} (ru_{\ell}) = \frac{\lambda_{\ell}}{r} u_{\ell} - rk^2\nu u_{\ell} - r\epsilon(r)  (\ell+1/2) \calF_{\ell} 
\end{equation}

We note from \eqref{eq:Fell} that $\calF_L$'s are the Fourier-Legendre coefficients of $U_N(r,t) F(|U_N(r,t)|^2)$. 
Clearly, there exist $\alpha,\beta$ such that  $|U_N|^2\in [\alpha,\beta]$. 
We assume $\alpha,\beta$ to be known.
Our strategy is to approximate $F$ using its Fourier-Chebyshev expansion:
\begin{equation}
\label{eq:Fcheb}
F(|U_N|^2)\approx \mathcal{P}_d\left(\frac{2|U_N|^2-\alpha-\beta}{\beta-\alpha}\right)=\sum_{k=0}^{2^d-1} a_kT_k\left(\frac{2|U_N|^2-\alpha-\beta}{\beta-\alpha}\right),
\end{equation}
where $T_k$'s are the Chebyshev polynomials defined by $T_k(\cos\phi)=\cos(k\phi)$. 
We need to evaluate the Fourier-Legendre coefficients $\calF_L$ of $U_NF(|U_N|^2)$ in terms of the $a_k$'s.
In Section~\ref{sec:legtocheb}, we describe a general procedure to accomplish this task efficiently.

Towards this goal, we note first using \eqref{eq:jacobiconv} that
$$
|U_N|^2(t) = \left( \sum_{\ell=0}^N u_{\ell} P_{\ell} (t) \right)^2 = \sum_{\ell=0}^{2N} d_L P_L(t), \qquad d_L=(\{u_\ell\}\star\{u_{\ell}\})(L),
$$
\begin{equation}
\label{eq:translegendre}
\frac{2|U_N(t)|^2-\alpha-\beta}{\beta-\alpha}=\frac{1}{\beta-\alpha}\left\{(2d_0-\alpha-\beta)P_0(t)+2\sum_{\ell=1}^{2N} d_L P_L(t)\right\}.
\end{equation}
Similarly, 
\begin{equation}\label{eq:cLPL} 
|U_N|^2 U_N  = \left(  \sum_{\ell=0}^{2N} d_\ell P_{\ell} (t)\right) \left( \sum_{\ell=0}^N u_{\ell'}  P_{\ell'}(t) \right) 
               = \sum_{L=0}^{3N} c_L P_L(t), \qquad c_L=(\{d_\ell\}\star\{u_\ell\})(L).
\end{equation}
By comparing \eqref{eq:cLPL} with \eqref{eq:Fell} we also have
$
  c_L = \frac{2L+1}{2} \calF_L.
$

We convert the system of second order ODEs \eqref{ODE_Fell} to first order ODEs as follows.

For $\ell = 0,\ldots,N$ let
$
    v_{\ell} = \frac{ d (r u_{\ell}) } {dr},
$
then the boundary conditions are
$
   v_{\ell}(R_0) = u_{\ell}(R_0) + R_0 \frac{d u_{\ell}}{dr}|_{r=R_0} = h_{\ell} + R_0 g_{\ell}.  
$

 Let
\[ 
\vec{Z} = [  Z_1 \;  Z_2 \; \cdots \; Z_{2N+2}]^\top =
   [ ru_0 \; ru_1 \; \cdots  \; ru_N \; v_0 \; v_1 \; \cdots \; v_N]^\top
\]
We can re-write the above system into the form $d\vec{Z}/dr = \mathfrak{F}(r,\vec{Z})$ with
\[
\mathfrak{F}(r,\vec{Z}) = 
\begin{bmatrix}
Z_{N+2} \\
Z_{N+3} \\
\vdots \\
Z_{2N+2} \\
\frac{\lambda_0}{r^2} Z_{1} 
 - k^2 \nu(r) Z_{1} - r\epsilon(r) (0+1/2) \calF_{0} \\  
\frac{\lambda_1}{r^2} Z_{2} 
 - k^2 \nu(r) Z_{2} - r\epsilon(r) (1+1/2) \calF_{1} \\  
\vdots \\
\frac{\lambda_L}{r^2} Z_{N+1} 
 - k^2 \nu(r) Z_{N+1} - r\epsilon(r) (L+1/2) \calF_{N}
\end{bmatrix}
\]
with the initial condition
\begin{align*}
Z(R_0) &= 
[ru_0(R_0) \; ru_1(R_0) \; \cdots \; r u_N(R_0) \; v_0 (R_0)\; v_1(R_0) \; \cdots \; v_N (R_0) ]  \\
&=
[R_0 h_0 \; R_0 h_1 \; \cdots \; R_0 h_N \;  h_0 + R_0 g_0 \;  h_1 + R_0 g_1 \; \cdots \; h_N + R_0 g_N]
\end{align*}

We may now use standard ODE solvers.
In our experiments we used the adaptive solver \textsf{ode45} in Matlab\textsuperscript{\textregistered}.

\section{The inverse problem}\label{sec:inverse}
In the inverse problem, we are given the values
$U(r_i)$ are known on the collection of points $\mathcal{R}:=\{r_i : i=1,\ldots, M\}$ which might not be equally spaced on the
interval $[R_0, R_1]$ since they might come from an adaptive ODE solver.  
The corresponding values $u_\ell(r_j)$ can be computed using numerical integration.
In our numerical experiments, we can extract $u_{\ell}$ directly from the numerical solutions of the ODE solver. 

Our approach is to evaluate $\calF_L$ first using \eqref{ode:uell}. 
In turn, this requires computing the second derivative of  $ru_{\ell}$ at $r=r_j$ for non-equidistant values $r_j$.
These are  computed  by 
\[
\left. \frac{\partial^2}{\partial^2 r} (ru_{\ell}(r) ) \right|_{r = r_j}
                \approx 
\frac{h^{-}_{j} r_{j+1}u_{\ell}(r_{j+1})  + h^{+}_{j} r_{j-1}u_{\ell}(r_{j-1}) - (h^{+}_j + h^{-}_j) r_j u_{\ell}(r_j)}            
                { 0.5 h^{-}_j h^{+}_j (h_j^{+} + h_j^{-})},
\]
with $h^{+}_j = r_{j+1} - r_{j}$ and $h^{-}_j = r_{j}-r_{j-1}$.
We then compute the approximated $\calF_{\ell}$ at $r = r_j$ via
\[
\calF_{\ell} = \frac{-2}{(2\ell+1)\epsilon(r)} \left( \frac{1}{r} \frac{\partial^2}{\partial^2 r} (r u_{\ell} ) 
 - \frac{\lambda_\ell}{r^2} u_{\ell}  + k^2 \nu(r)  u_{\ell} \right)
\]
The next task is to approximate $F$ from the values of $\calF_L$'s. 
Since we know $\calF_L$'s, this leads to a  (not necessarily square) system of non-linear equations. 
In turn, the $a_k$'s are determined using a least squares computation.
Thus, the problem reduces to computing $a_k$'s using the expansion \eqref{eq:translegendre}.

\section{Numerical experiments}\label{sec:numerics}
The expansion of a plane wave is given by Morse and Ingard \cite{MorIng68}
\begin{equation}
e^{i \bk \cdot \br} = 
\sum_{\ell=0}^\infty (2\ell+1) i^{\ell} 
 P_\ell(\widehat{\bk} \cdot \widehat{\br}) j_{\ell}(kr),
\end{equation}
where $\widehat{\bk} = \bk / \|\bk\|$,
$\widehat{\br} =  \br /\|\br\|$, $P_{\ell}(t)$ is
the Legendre polynomial of degree $\ell$ and $j_{\ell}(kr)$ is the $\ell$th
spherical Bessel function of the first kind.  
Here $\br$ is the position vector of length $r$, $\bk$ is the wave
vector of length $k$. In the special case when $\bk$ is aligned with
the $z$-axis, we have
\[
 e^{ikr \cos\theta} = \sum_{\ell=0}^\infty (2\ell+1) i^{\ell} 
 P_\ell(\cos\theta) j_{\ell}(kr),
\]
where $\theta$ is the spherical polar angle of $\br$. With $t=\cos\theta$, we have $H = e^{ikR_0 t}$ and
\[
  h_{\ell} = (2\ell+1) i^{\ell} j_{\ell}(k R_0)
\]
and by using the identity $\frac{d}{dz}j_{\ell}(z) = j_{\ell-1}(z) - \frac{(\ell+1)}{z} j_{\ell}(z)$, we have
\begin{align*}
  g_{\ell} &= (2\ell+1) i^{\ell} \left. \frac{\partial j_{\ell}(kr) }{\partial r} \right|_{r=R_0} 
   =(2\ell+1) i^{\ell} 
  \frac{1}{k}\left( j_{\ell-1}(kR_0) - \frac{\ell+1}{kR_0} j_{\ell}(kR_0)\right) .
\end{align*}

%
\subsection{Experiment 1}
\begin{figure}[h]
\centering
\begin{tabular}{cc}
\includegraphics[width=0.5\textwidth]{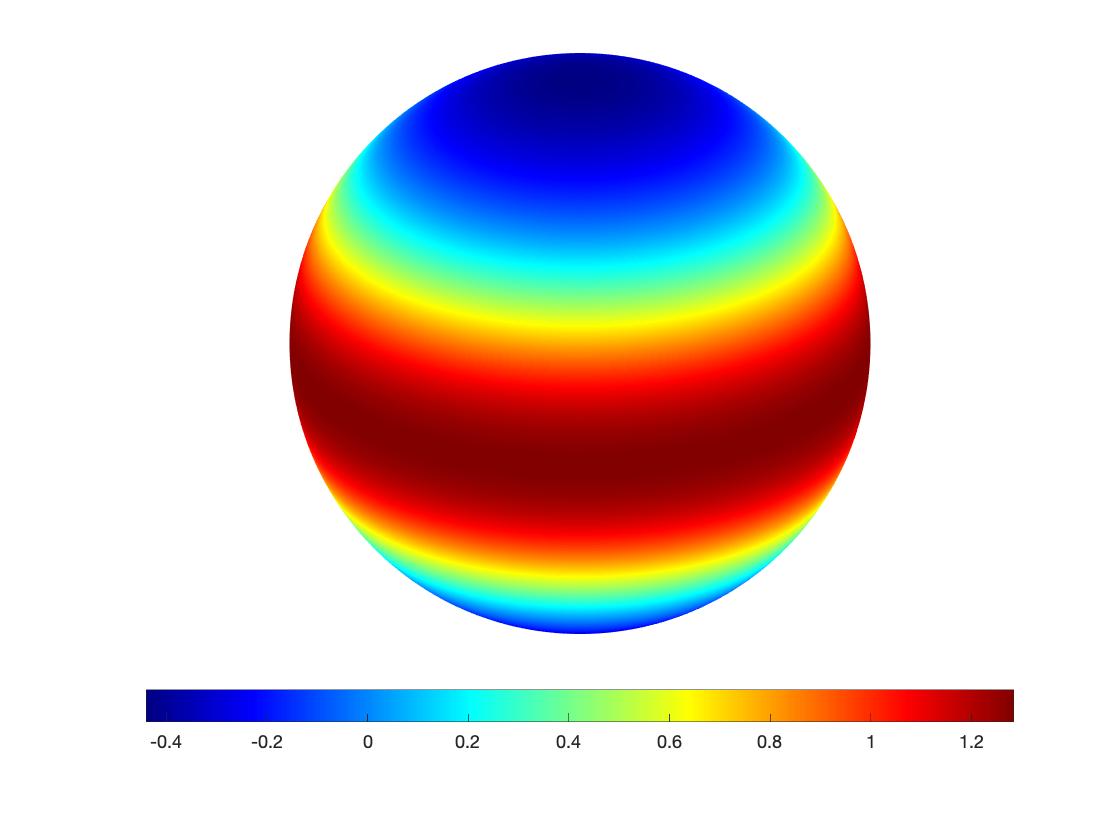} & 
\includegraphics[width=0.5\textwidth]{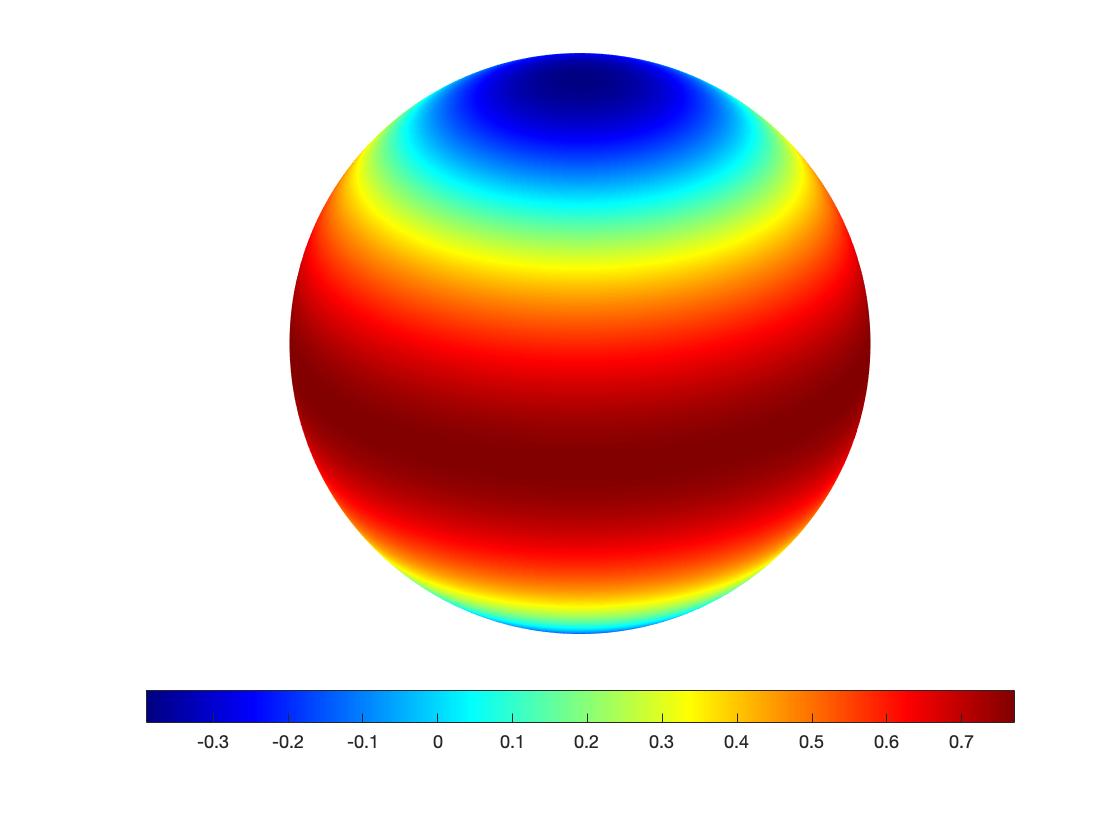}
\end{tabular}
\caption{$U(r=R_1)$ with $R_1=2$ for $\epsilon=2$, $\nu=0.1$ and $k=1$ for Experiment 1 (left panel) and Experiment
2 (right panel).}\label{fig:U80 sol experiment1}
\end{figure}

We consider the forward problem
\begin{equation}
\Delta U(\bx) + k^2 \nu U(\bx) = -\epsilon |U(\bx)|^4 U(\bx),
\end{equation}
where $k,\nu$ and $\epsilon$ are positive constants
on the spherical shell $\Omega$ with inner radius $R_0=1$ and outer radius $R_1=2$.
The boundary conditions on the inner sphere are given by
\[
 U(R_0) = e^{ikR_0 t}, \quad \frac{\partial U}{\partial r} = \frac{\partial}{\partial r} e^{ikrt}|_{r=R_0}, \quad t = \cos\theta
\] 
The numerical solution $U(R_1)$ of the forward problem is given in the left panel of Figure~\ref{fig:U80 sol experiment1}.

We now consider the inverse problem. On the right-hand side, in our framework $F(|U|^2) = |U|^4$, so $F(t) = t^2$. The function $F$
can be expressed as a linear combination of Chebyshev polynomials $T_0$ and $T_2$ as
\[
 F(t) = \frac{1}{2} T_0(t) + \frac{1}{2} T_2 (t). 
\] 
So the exact coefficients are $a_0 = 0.5$, $a_1 = 0$ and $a_2 = 0.5$. 
The computed coefficients from the inverse problem $a_0, a_1, a_2$ on each ring are shown in Table~\ref{tab:FU4coeff}.

\begin{table}[h]
\centering
\begin{tabular}{|c|c|c|c|}
\hline
    $r$  & $a_0$ & $a_1$  & $a_2$ \\
    \hline
1.0009 &  5.0000e-01 & -6.3171e-06 &  5.0000e-01\\
1.0018 &  5.0000e-01 &-6.2524e-06 &  5.0000e-01\\
1.0027 &  5.0000e-01 & -6.3744e-06 &   5.0000e-01\\
1.0036 &  4.9591e-01 & -7.1347e-03 &  4.9815e-01\\
1.0065 &  4.9995e-01 &  -6.3584e-05&   4.9997e-01\\
1.0094 &  4.9996e-01 & -6.1446e-05 &  4.9997e-01\\
1.0123 &  4.9995e-01 & -6.5269e-05 &  4.9997e-01\\
1.0152 &  4.9992e-01 & -1.2534e-04  & 4.9995e-01\\
1.0181 &  4.9995e-01 & -6.4420e-05  & 4.9997e-01\\
\hline
\end{tabular}
\caption{Computed Chebyshev coefficients for $F(|U|^2) = |U|^4$} \label{tab:FU4coeff}
\end{table}

\subsection{Experiment 2}
Let $F(|U|^2) = \sin |U|^2$ and $|U|^2 \in [\alpha, \beta]$.
Let
\[
  q_k (z) := \begin{cases}
            J_0 (z), & \text{ if } k = 0,\\
            2 J_k(z), &\text{ if } k \in \N,
            \end{cases}
\]
where $J_k$ is the Bessel's function of order $k$.
From Watson's book, \cite[page 22,(3)-(4)]{Watson}, we have for $t \in [-1,1]$,
\begin{align*}
\sin (\gamma + zt) &= \sin \gamma \cos(zt) + \cos \gamma \sin(zt) \\
    &= \sin\gamma \sum_{k=0}^\infty (-1)^k q_{2k}(z) T_{2k}(t) + \cos\gamma  \sum_{k=0}^\infty (-1)^k q_{2k+1}(z) T_{2k+1}(t) \\
    &=  \sum_{k=0}^\infty \sin \left(\gamma+ \frac{2k \pi}{2} \right) q_{2k}(z) T_{2k}(t) 
    + \sum_{k=0}^\infty \sin \left(\gamma+ \frac{(2k +1)\pi}{2} \right)  q_{2k+1}(z) T_{2k+1}(t) \\
     &= \sum_{n=0}^\infty \sin \left( \gamma + \frac{n\pi}{2} \right) q_n(z) T_n(t).
\end{align*}
So with $\gamma = (\alpha+\beta)/2$ and $z=(\beta-\alpha)/2$, $|U|^2 = \gamma+ zt$, then $t \in [-1,1]$ and
\[
\sin(|U|^2) = \sum_{n=0}^\infty \sin \left( \frac{\alpha+\beta}{2} + \frac{n\pi}{2}\right) q_n \left(  \frac{\beta-\alpha}{2}\right) T_n(t),
\quad t = \frac{2|U|^2-\alpha-\beta}{\beta-\alpha}.
\]

Let's assume $|U| \in [0,1]$, that is, $\alpha = 0$, $\beta = 1$ and we use only the first $8$ terms of the infinite series above to define
\begin{equation}\label{eq:FsinU}
F(|U|^2) = \sum_{n=0}^7 \sin \left( \frac{1}{2} + \frac{n\pi}{2}\right) q_n \left(  \frac{1}{2}\right) T_n(t),
        \quad t = 2|U|^2-1.
\end{equation}
So the exact coefficients are $a_n = \sin(1/2+n \pi/2) q_n(1/2)$ for $n=0,\ldots,7$. 
The numerical solution of the forward problem $U(R_1)$ is given in right panel of Figure~\ref{fig:U80 sol experiment1}.
For the inverse problem, the computed
coefficients $a_n$ for $n=0,\ldots,7$ on each ring are shown in Table~\ref{tab:FsinUcoeff}.


\begin{table}[h]
\centering
\begin{tabular}{|c|c|c|c|c|}
\hline
        & $a_0$      &  $a_1$      &  $a_2$      & $a_3$ \\
exact   & 4.4993e-01 &  4.2522e-01 & -2.9345e-02 &   -4.4998e-03 \\
\hline
$r=$1.001634 & 4.4993e-01 &  4.2522e-01 & -2.9344e-02 &   -4.4999e-03\\
$r=$1.003268 & 4.4993e-01 &  4.2522e-01 & -2.9344e-02 &   -4.4999e-03\\
$r=$1.004902 & 4.4993e-01 &  4.2522e-01 & -2.9344e-02 &   -4.4999e-03\\
\hline
\end{tabular}
\begin{tabular}{|c|c|c|c|c|}
\hline
         & $a_4$ & $a_5$ & $a_6$ & $a_7$ \\
exact    & 1.5412e-04 &  1.4135e-05 & -3.2224e-07 &   -2.1090e-08 \\
\hline
$r=$1.001634 & 1.5409e-04 & 1.4148e-05 & -3.2558e-07 &  -2.0502e-08\\
$r=$1.003268 & 1.5408e-04 & 1.4150e-05 & -3.2638e-07 &  -2.0376e-08\\
$r=$1.004902 & 1.5408e-04 & 1.4150e-05 & -3.2587e-07 &  -2.0499e-08\\
\hline
\end{tabular}
\caption{Computed Chebyshev coefficients for $F(|U|^2)$ as in \eqref{eq:FsinU}} \label{tab:FsinUcoeff}
\end{table}



\section{Computational issues}\label{sec:legtocheb}
Let $f \in C([-1,1]$, 
$$
\mathcal{P}_d(t)=\sum_{k=0}^{2^d-1} a_kT_k(t), \quad \mathbb{P}_d(t)=\mathcal{P}_d(f(t)), \qquad t\in [-1,1].
$$
We wish to compute the Fourier-Legendre coefficients $\{b_\ell\}$ of $\mathbb{P}_d$ explicitly and efficiently using the Fourier-Legendre coefficients of $f$ and the coefficients $a_k$.

We proceed inductively. 
If $d=1$, then we observe that
\begin{equation}
\label{eq:initialization}
\begin{aligned}
\langle T_0(f), P_0\rangle&=1, \qquad \langle T_1(f), P_1\rangle=\hat{f}(1)\\
\mathbb{P}_1(t)&=\frac{1}{2}+\frac{3}{2}\hat{f}(1)P_1(t).
\end{aligned}
\end{equation}

Next, we assume that the problem is solved in the case of polynomials of degree $\le 2^{d-1}-1$. 
Using the recurrence relations
\begin{equation}\label{eq:chebrecur}
 T_{2^j+k}= 2T_{2^j}T_k -T_{2^j-k}, \qquad j=1,2,\cdots, \ k=1,\cdots 2^j,
\end{equation}
it is not difficult to deduce that
\begin{equation}
\label{eq:dyadicrec}
\begin{aligned}
\sum_{k=0}^{2^d-1}a_kT_k&=\sum_{k=0}^{2^{d-1}-1}(a_k-a_{2^d-k})T_k +2T_{2^{d-1}}\sum_{k=0}^{2^{d-1}}a_{k+2^{d-1}}T_k\\
&=\mathcal{Q}_{d-1}+2T_{2^{d-1}}\widetilde{\mathcal{R}}_{d-1}
\end{aligned}
\end{equation}
for polynomials $\mathcal{Q}_{d-1}$, $\widetilde{\mathcal{R}}_{d-1}$ of degree at most $2^{d-1}$. 
We let $\mathbb{Q}_{d-1}(t)=\mathcal{Q}_{d-1}(f(t))$ and $\widetilde{\mathbb{R}}_{d-1}(t)=\widetilde{\mathcal{R}}_{d-1}(f(t))$.
Given our induction hypothesis, we may now compute
\begin{equation}
\label{eq:finalsol}
\widehat{\mathbb{P}_d}= \widehat{\mathbb{Q}_{d-1}}
+2(\widehat{T_{2^{d-1}}\circ f})\star\widehat{\widetilde{\mathbb{R}}_{d-1}}.
\end{equation}
Using \eqref{eq:initialization} and \eqref{eq:finalsol}
one can compute $\widehat{\mathbb{P}_d} = \widehat{F}$ with $O(d)$ convolutions.

\paragraph{Acknowledgements}
The authors thank the support of the Australian Research Council,  
Q.L.G. was supported by DP180100506.
The research of HNM was supported in part by ARO grant W911NF2110218 and NSF DMS grant 2012355.
\ifx\printbibliography\undefined
    \bibliographystyle{plain}
    \bibliography{nonlinHelm}

\begin{thebibliography}{1}

\bibitem{askeyspecial}
R.~Askey.
\newblock {\em Orthogonal polynomials and special functions}.
\newblock SIAM, 1975.

\bibitem{BarFibTsy07}
G.~Baruch, G.~Fibich, and S.~Tsynkov.
\newblock High-order numerical method for the nonlinear {Helmholtz} equation
  with material discontinuities in one space dimension.
\newblock 2007.

\bibitem{FibTsy01}
G.~Fibich and S.~Tsynkov.
\newblock High-order two-way artificial boundary conditions for nonlinear wave
  propagation with backscattering.
\newblock {\em Journal of Computational Physics}, 171:632--677, 2001.

\bibitem{FibTsy05}
G.~Fibich and S.~Tsynkov.
\newblock Numerical solution of the nonlinear helmholtz equation using
  nonorthogonal expansions.
\newblock {\em Journal of Computational Physics}, 210:183--224, 2005.

\bibitem{MorIng68}
P.~M. Morse and K.~U. Ingard.
\newblock {\em Theoretical Acoustics Vols. 1–2}.
\newblock McGraw-Hill Book Company, 1968.

\bibitem{Watson}
G.~N. Watson.
\newblock {\em A treatise on the theory of Bessel functions}.
\newblock Cambridge Mathematical Library, 1996.

\end{thebibliography}
\else\printbibliography\fi
\end{document}